\documentclass[11pt]{article}
\usepackage[T1]{fontenc}
\usepackage[utf8]{inputenc}
\usepackage{lmodern}
\usepackage{xcolor}
\usepackage{amsmath,amssymb,amsthm,mathtools}
\usepackage[margin=1in,footskip=30pt]{geometry}
\usepackage{microtype}
\usepackage[hidelinks]{hyperref}
\newtheorem{theorem}{Theorem}[section]
\newtheorem{lemma}[theorem]{Lemma}

\newtheorem{corollary}[theorem]{Corollary}
\newtheorem{conjecture}[theorem]{Conjecture}
\theoremstyle{definition}
\newtheorem{definition}[theorem]{Definition}
\theoremstyle{remark}

\newcommand{\E}{\mathbb E}
\newcommand{\Prob}{\mathbb P}
\newcommand{\supp}{\operatorname{supp}}
\newcommand{\claimqed}{%
  {\renewcommand{\qedsymbol}{$\blacksquare$}\qed}%
}

\title{Ore-degree conditions for hypergraph matchings}
\author{Mingyang Guo\thanks{School of Mathematics and Statistics, Zhengzhou University, Zhengzhou, Henan, China}, Tianhong Kan\thanks{School of Mathematics and Statistics, Xi’an Jiaotong University, Xi’an, Shaanxi, China.}, Hongliang Lu\footnotemark[2] and Feihong Yuan\footnotemark[2]}

\date{}
\hypersetup{pdftitle={Ore-degree conditions for hypergraph matchings},pdfauthor={Anonymous},pdfsubject={Hypergraph matchings and Ore-degree conditions}}
\begin{document}
\maketitle
\begin{abstract}
Let $\sigma_k(H)$ be the minimum sum of the vertex degrees over all non-edges of an $n$-vertex $k$-uniform hypergraph $H$. For every $k\ge3$, $n\geq ks$ and sufficiently large $n$, we establish a sharp Ore-type degree condition that forces a matching of size $s$ in $H$. Specifically, we prove that if 
\[
 \sigma_k(H)>\max\left\{
 k\left(\binom{n-1}{k-1}-\binom{n-s}{k-1}\right),
 (k-c)\binom{ks-2}{k-1}\right\},
\]
then $H$ has a matching of size $s$, where $c=n-ks+1$. In particular, this proves a conjecture of Balogh, Palmer and Raeisi for $k\geq 3$ and sufficiently large $n$.
\end{abstract}

\section{Introduction}
A \emph{hypergraph} $H$ is a pair $H=(V,E)$, where $V:=V(H)$ is a set of vertices and $E:=E(H)$ is a set of non-empty subsets of $V$. For a finite set $S$ and a positive integer $k$, write $\binom Sk=\{X\subseteq S:|X|=k\}$. When $E(H)\subseteq\binom{V(H)}k$, we call $H$ a \emph{$k$-uniform hypergraph}, or a \emph{$k$-graph}. Write $e(H):=|E(H)|$. For $v\in V(H)$, we denote by $d_H(v)$ the number of edges containing $v$ in $H$. Let $\delta_1(H):=\min\{d_H(v):v\in V(H)\}$ denote the \emph{minimum vertex degree} of $H$.  

A \emph{matching} in $H$ is a set of pairwise disjoint edges. We use $\nu(H)$ to denote the size of a maximum matching in $H$. A matching is \emph{perfect} if it covers $V(H)$. For $n\geq ks$, an $n$-vertex complete $k$-graphs has a matching of size $s$, so we only consider noncomplete $k$-graph in this paper. We define the \emph{Ore-degree} of a  $k$-graph $H$ by
 $\sigma_k(H)=\min\left\{\sum_{v\in e}d_H(v):
 e\in\binom Vk\setminus E(H)\right\}$.

In 1965, Erd\H{o}s~\cite{Erdos}  asked for the determination of the maximum possible number of edges that can appear in any $k$-graph $H$ with $\nu(H)\leq s$ and proposed the celebrated Erd\H os matching conjecture.
\begin{conjecture}
    [Erd\H os Matching Conjecture, \cite{Erdos}]\label{EMCconj}
	Let $n,s,k$ be three positive integers such that $k\geq2$ and $n\geq k(s+1)-1$. If $H$ is a $k$-graph on $[n]$ and $\nu(H)\leq s$, then
	\begin{equation*}
	e(H)\leq\max\left\{\binom{n}{k}-\binom{n-s}{k},\binom{k(s+1)-1}{k}\right\}.
	\end{equation*}
\end{conjecture}
For progress on the Erd\H{o}s Matching Conjecture, see, e.g.,
\cite{Frankl2013,Frankl2017a,Frankl2017b,FK,KK2023,CLZ}.

Bollob\'as, Daykin and Erd\H{o}s~\cite{BDE} considered a minimum vertex degree version of Conjecture \ref{EMCconj} and proved that if
 $\delta_1(H)>\binom{n-1}{k-1}-\binom{n-s}{k-1}$,
then $\nu(H)\ge s$
when $n>2k^3(s+1)$. Huang and Zhao~\cite{HZ} extended the range to $n\ge3k^2s$. Guo, Lu and Jiang~\cite{GLJ}  proved that, for every integer
$k\ge3$ and every $\beta>0$, the same conclusion holds   whenever $1\le s\le(\frac{k}{2(k-1)}-\beta)\frac nk$  and $n$ is sufficiently large. For $k=3$, K\"uhn, Osthus and Treglown~\cite{KOT} obtained the sharp bound for the existence of a matching of size $s$ whenever $1\le s\le n/3$ and $n$ is sufficiently large.

For perfect matchings, H\`an, Person and Schacht~\cite{HPS} obtained an asymptotically sharp minimum vertex degree bound for $3$-graphs, and Khan~\cite{Khan3,Khan4} determined the exact thresholds for $k=3,4$. Recently, Han, Lu, Wang and Yuan~\cite{HLWY} determined the sharp minimum vertex degree threshold for every $k\ge3$ and sufficiently large $n$, as a special case of their minimum $d$-degree theorem. For further results on degree conditions and fractional matchings, we refer the reader to~\cite{AFHRRS,FK,Han2016,KOTown,LYY,TZ2016}.

For graphs of even order, the classical Dirac's theorem~\cite{Dirac}  and Ore's theorem~\cite{Ore} give sufficient degree conditions for perfect matchings. Dirac's theorem implies that every graph $G$ of even order $n\ge4$ with $\delta_1(G)\ge n/2$ has a perfect matching. Ore's theorem gives the same conclusion under the weaker assumption $\sigma_2(G)\ge n$. This motivates replacing a minimum vertex degree condition by a condition on degree sums over non-edges.
Recently, Balogh, Palmer and Raeisi~\cite{BPR} considered the Ore-degree condition over non-edges and proved the following theorem.
\begin{theorem}[{Balogh--Palmer--Raeisi, \cite[Theorem 1.5]{BPR}}]\label{thm:bpr}
	Let $n,s,k$ be integers such that  $k\ge3$, $s\ge2$, and $n\ge3k^2(s-1)$. If an $n$-vertex $k$-graph $H$ satisfies $\sigma_k(H)>k\left(\binom{n-1}{k-1}-\binom{n-s}{k-1}\right)$, then $\nu(H)\ge s$.
\end{theorem}
The Ore-degree bound is best possible. To see it, we define the $k$-graph $H_{k}(U,W)=(V,E)$, where $V$ is partitioned into $U\cup W$, and  $E=\{e\in \binom{V}{k}\ :\ 1\leq |e\cap W|\leq k-1\}$.
When $|W|=s-1$, we denote $H_{k}(U,W)$ by $H_k(n,s)$. Let $D(n,k,s)=\binom{n-1}{k-1}-\binom{n-s}{k-1}$. Since every edge intersects $W$, $\nu(H_{k}(n,s))=s-1$.
For $u\in U$ and $w\in W$,
$ d_{H_k(n,s)}(w)>d_{H_k(n,s)}(u)
 =D(n,k,s)$.
Thus $\sigma_k(H_{k}(n,s))=kD(n,k,s)$.

They proposed the following conjecture. 

\begin{conjecture}[{Balogh--Palmer--Raeisi, \cite[Conjecture 2]{BPR}}]\label{conj:BPR}
Let $n,s,k$ be integers such that $k\ge3$, $s\ge2$ and $n>2ks$. Every $n$-vertex $k$-graph $H$ satisfying
\[
 \sigma_k(H)>k\left(\binom{n-1}{k-1}-\binom{n-s}{k-1}\right)
\]
contains a matching of size $s$.
\end{conjecture}
Subsequently, Zhang, Cao and Lu \cite{ZCL} proved Conjecture \ref{conj:BPR} for $n\geq 3ks$ and sufficiently large $s$.
In this paper, we determine the sharp Ore-degree bound for all $n\ge ks$ and sufficiently large $n$. 
\begin{theorem}\label{thm:main}
For every integer $k\ge3$, there exists $n_0=n_0(k)$ such that the following holds for all integers $n\ge n_0$. Let $s,c$ be integers such that $s\geq 1$, $n\geq ks$ and $c=n-ks+1$. Suppose that $H$ is an $n$-vertex $k$-graph. If
\[
 \sigma_k(H)>\max\left\{
 k\left(\binom{n-1}{k-1}-\binom{n-s}{k-1}\right),
 (k-c)\binom{ks-2}{k-1}\right\},
\]
then $H$ contains a matching of size $s$. 
\end{theorem}

The bound in Theorem~\ref{thm:main} is best possible.
The $k$-graph $H_k(n,s)$ has matching number $s-1$ and
Ore-degree $kD(n,k,s)$.
Let $B_k(n,s)$ consist of a complete $k$-graph on $ks-1$
vertices together with $c=n-ks+1$ isolated vertices.
Then $\nu(B_k(n,s))=s-1$ and
$\sigma_k(B_k(n,s))
=\max\{k-c,0\}\binom{ks-2}{k-1}$.
Since $D(n,k,s)\ge0$, the larger of the Ore-degrees of these
two constructions is exactly the bound in
Theorem~\ref{thm:main}.


By Theorem \ref{thm:main}, Conjecture \ref{conj:BPR} holds for all sufficiently large $n$.

\begin{corollary}\label{cor:BPR}
	For every $k\ge3$, Conjecture~\ref{conj:BPR} holds for all sufficiently large $n$.
\end{corollary}

\medskip
\noindent\textbf{Proof strategy.}
We first consider the case $\frac{n-c_0}{k}\leq s\leq \frac{n}{k}$ for some fixed constant $c_0$ and assume
that $\delta_1(H)\ge\alpha\binom{n-1}{k-1}$ for some
fixed $\alpha>0$.
After deleting $c-1$ vertices, we distinguish two cases, where $c=n-ks+1$.
If the remaining hypergraph contains a large sparse
vertex set in which most vertices have degrees close
to $D(ks,k,s)$, we exploit this structure to construct
an $s$-matching in $H$
(Section~\ref{sec:extremal}).
Otherwise, a compactness argument using approximate
fractional covers and the equality case of a probability
inequality yields an almost perfect matching
(Section~\ref{sec:endpoint}).
We combine this with the absorbing lemma
(Section~\ref{sec:absorb}) to obtain a perfect matching
(Section~\ref{sec:non-extremal}).

We then remove the minimum degree assumption by finding
a small matching covering all but at most $c-1$ low-degree
vertices. After deleting this matching and a suitable
set of $c-1$ vertices, we apply the preceding result
to the remaining hypergraph.
The clique term in the Ore-degree threshold is crucial
in this step (Section~\ref{sec:extension}).

Finally, we extend the result to $1\le s\le n/k$ by adjoining
$a$ new vertices and all $k$-sets intersecting them, with $a$ chosen
so that $(N-c_0)/k\le S\le N/k$, where $N=n+a$ and $S=s+a$. Applying the preceding result to the enlarged hypergraph
and deleting the matching edges that intersect the new vertices
then yields an $s$-matching in $H$
(Section~\ref{sec:extension}).


\smallskip
We finish this section with notation used below. Let $[n]=\{1,\ldots,n\}$ and $\binom ab=0$ for integers $a,b$ with $b<0$ or $a<b$. For $X\subseteq V(H)$, the induced subhypergraph $H[X]$ has vertex set $X$ and edge set $\{e\in E(H):e\subseteq X\}$. Write $H-X=H[V(H)\setminus X]$ . For a vertex $v$, its link is
$ N_H(v)=\left\{S\in\binom{V(H)\setminus\{v\}}{k-1}:S\cup\{v\}\in E(H)\right\}$.
For a matching $M$, write $V(M)=\bigcup_{e\in M}e$.

\section{Preliminaries}\label{sec:prelim}

We need the following two results proved by Cao, Liu and Zhang~\cite{CLZ} to find an almost perfect matching in the non-extremal case.
\begin{lemma}[{Cao--Liu--Zhang, \cite[Lemma 3.2]{CLZ}}]\label{lem:approx}
Fix $k\ge2$, $0\le\xi<1/k$, and $\rho>0$. There exists $n_0$ such that, for every $n\ge n_0$, the following holds. If $0\le x\le\xi$ and $H$ is an  $n$-vertex $k$-graph satisfying $\nu(H)\le xn$, then there is a function $w:V(H)\to[0,1]$ such that
\[
 \sum_{v\in V(H)}w(v)\le(x+\rho)n
\]
and $\sum_{v\in e}w(v)\ge1$ for all but at most $\rho n^k$ edges $e\in E(H)$.
\end{lemma}

\begin{theorem}[{Cao--Liu--Zhang~\cite[Theorem 2.5]{CLZ}}]\label{thm:tail}
Let $r\ge1$ and $T\ge r+1$. If $Y_1,\ldots,Y_r$ are independent nonnegative random variables with $\E Y_i\le1$, then
\[
 \Prob(Y_1+\cdots+Y_r\ge T)\le1-\left(1-\frac1T\right)^r.
\]
Equality holds if and only if each $Y_i$ takes values $0,T$ with probabilities $1-1/T,1/T$, respectively.
\end{theorem}

We shall also use the following special case of a theorem of H\`an, Person and Schacht \cite{HPS} to deal with the extremal case. 
\begin{theorem}[{H\`an--Person--Schacht, \cite[Theorem 2]{HPS}}]\label{thm:hps}
	Let $J$ be a $k$-uniform, $k$-partite hypergraph whose $k$ partition classes each have size $n$. If
	$\delta_1(J)>\frac{k-1}{k}n^{k-1}+kn^{k-2}$,
	then $J$ has a perfect matching.
\end{theorem}

The following special case of the absorption lemma of Lo and Markstr\"om allows us to construct an absorbing set for vertices uncovered by the almost perfect matching.
\begin{lemma}[{Lo--Markstr\"om, \cite[Lemma 1.1]{LoMarkstrom}}]\label{lem:LMabsorb}
	Let $k\ge2$ and $i\ge1$ be integers, and let $\beta>0$. There is an integer $n_0=n_0(k,i,\beta)$ such that the following holds for $n\ge n_0$. Suppose that $H$ is a $k$-graph on $n$ vertices and, for every two distinct vertices $x,y\in V(H)$, there are at least $\beta n^{ik-1}$ sets $X\subseteq V(H)\setminus\{x,y\}$ of size $ik-1$ for which both $H[X\cup\{x\}]$ and $H[X\cup\{y\}]$ have perfect matchings. Then there exists $A\subseteq V(H)$ of size $|A|\le\frac{(\beta/2)^k}{4ik(k-1)}\,n$ with $|A|\in k\mathbb Z$ such that $H[A\cup W]$ has a perfect matching for every $W\subseteq V(H)\setminus A$ of size
	$|W|\le\frac{(\beta/2)^{2k}}{32i^2k(k-1)^2}\,n$ with $|W|\in k\mathbb Z$.
\end{lemma}

\section{An absorbing lemma}\label{sec:absorb}
Using Lemma~\ref{lem:LMabsorb}, we obtain the following absorbing lemma under the minimum-degree and the Ore-degree conditions. 

\begin{lemma}\label{lem:absorb}
	Let $k$ be an integer such that $k\ge3$. Let $\alpha,\zeta$ be constants such that $\alpha,\zeta\in(0,1)$, and let $\beta=(\alpha\zeta/(512k^2))^{2k-1}/(k^2-1)!$. There exists $n_0=n_0(k,\alpha,\zeta)$ such that the following holds for every $n\ge n_0$. Suppose that $H$ is an $n$-vertex $k$-graph with $\delta_1(H)\ge\alpha\binom{n-1}{k-1}$ and $\sigma_k(H)\ge(k/2+\zeta)\binom{n-1}{k-1}$. Then there exists $A\subseteq V(H)$ of size $|A|\le(\beta/2)^kn/(8k(k-1))$ such that $H[A\cup W]$ has a perfect matching for every $W\subseteq V(H)\setminus A$ satisfying $|W|\le(\beta/2)^{2k}n/(32k^3(k-1)^2)$ and $k\mid |W|$.
\end{lemma}

\begin{proof}
	Let $h=\zeta/(8k)$ and $\tau=\alpha\zeta/(128k^2)$. Thus $\beta=(\tau/4)^{2k-1}/(k^2-1)!$.  Let $L=\{x\in V(H):d_H(x)<(1/2+h)\binom{n-1}{k-1}\}$. We distinguish two cases.
	
	\smallskip\noindent
	\emph{Case 1. $|L|\geq\tau n$.}
	
	Since $\sum_{x\in S}d_H(x)\leq(k/2+kh)\binom{n-1}{k-1}<\sigma_k(H)$ for every $S\in\binom{L}{k}$, $H[L]$ is a complete $k$-graph. Let $L_1$ be a subset of $L$ such that $|L_1|=\lceil |L|/2\rceil$ and $d_H(x)\geq d_H(y)$ for every $x\in L_1$ and every $y\in L\setminus L_1$. Fix distinct vertices $v,w\in V(H)$. Choose $u\in L_1\setminus\{v,w\}$.
	
	Let $x\in\{v,w\}$ and $L_2=L\setminus L_1$. If $\binom{L_2\setminus\{x\}}{k-1}\subseteq N_H(x)$, we have $$|N_H(u)\cap N_H(x)|\geq\left|\binom{L_2\setminus\{x\}}{k-1}\right|\geq(\tau n/3)^{k-1}/(k-1)!,$$ since $u\notin L_2$, $H[L]$ is complete, and $|L_2|\geq\tau n/2-1$. Otherwise, there exists $S\in\binom{L_2\setminus\{x\}}{k-1}$ such that $S\cup\{x\}\notin E(H)$, we have $d_H(x)+\sum_{y\in S}d_H(y)\geq\sigma_k(H)$. Thus $d_H(x)\geq\sigma_k(H)-\sum_{y\in S}d_H(y)\geq\sigma_k(H)-(k-1)d_H(u)$ and $$|N_H(u)\cap N_H(x)|\geq d_H(u)+d_H(x)-\binom n{k-1}\geq\sigma_k(H)-(k-2)d_H(u)-\binom n{k-1}\geq h\binom{n-1}{k-1},$$ where the last inequality follows from $d_H(u)<(1/2+h)\binom{n-1}{k-1}$, $h=\zeta/(8k)$ and the fact that $n$ is sufficiently large. Since $(\tau/3)^{k-1}<h/2$, it follows that in either case, $|N_H(u)\cap N_H(x)|\geq(\tau n/3)^{k-1}/(k-1)!$.
	
	Choose disjoint sets $S\in N_H(v)\cap N_H(u)$ and $S'\in N_H(w)\cap N_H(u)$ avoiding $u,v,w$. Let $X=\{u\}\cup S\cup S'$. Thus $|X|=2k-1$, and the matchings $\{S\cup\{v\},S'\cup\{u\}\}$ and $\{S\cup\{u\},S'\cup\{w\}\}$ cover $X\cup\{v\}$ and $X\cup\{w\}$, respectively. There are at least $|L_1|-2\geq\tau n/4$ choices for $u$. For $S$ and $S'$, the number of choices at each step is at least $(\tau n/3)^{k-1}/(k-1)!-(2k-1)n^{k-2}\geq(\tau n/4)^{k-1}/(k-1)!$. It follows that there are at least $(\tau n/4)^{2k-1}/((k-1)!)^2$ choices of $(u,S,S')$. For each fixed $X\subseteq V(H)$ of size $2k-1$, the number of distinct triples $(u,S,S')$ satisfying $X=\{u\}\cup S\cup S'$ is at most $(2k-1)\binom{2k-2}{k-1}=(2k-1)!/((k-1)!)^2$. Therefore, every pair of distinct vertices $v,w$ has at least $(\tau/4)^{2k-1}n^{2k-1}/(2k-1)!\geq\beta n^{2k-1}$ sets $X\subseteq V(H)\setminus\{v,w\}$ of size $2k-1$ such that both $H[X\cup\{v\}]$ and $H[X\cup\{w\}]$ have perfect matchings.
	
		Apply Lemma~\ref{lem:LMabsorb} with $i=2$ and $\beta$. Then there exists $A\subseteq V(H)$ with $|A|\leq(\beta/2)^kn/(8k(k-1))$ and $|A|\in k\mathbb{N}$ such that $H[A\cup W]$ has a perfect matching for every $W\subseteq V(H)\setminus A$ satisfying $|W|\leq(\beta/2)^{2k}n/(128k(k-1)^2)$ and $|W|\in k\mathbb{N}$. This completes Case~1.
	
	\smallskip\noindent
	\emph{Case 2. $|L|<\tau n$.}
	
	Let $Q=V(H)\setminus L$. Thus $|Q|>(1-\tau)n$. For every $v\in V(H)$, we have $$\left|N_H(v)\cap\binom Q{k-1}\right|\geq\alpha\binom{n-1}{k-1}-\tau n\binom{n-2}{k-2}\geq(\alpha/2)\binom{n-1}{k-1}.$$ For distinct $u,z\in Q$, $$\left|N_H(u)\cap N_H(z)\cap\binom Q{k-1}\right|\geq2h\binom{n-1}{k-1}-\binom{n-1}{k-2}-\tau n\binom{n-2}{k-2}\geq h\binom{n-1}{k-1}.$$ 
	
	Fix distinct vertices $v,w$. There are at least $(\alpha/2)\binom{n-1}{k-1}-\binom{n-2}{k-2}$ choices for $S\in N_H(v)\cap\binom Q{k-1}$ avoiding $v,w$, and then at least $(\alpha/2)\binom{n-1}{k-1}-k\binom{n-2}{k-2}$ choices for $S'\in N_H(w)\cap\binom Q{k-1}$ avoiding $S\cup\{v,w\}$. 
	
	For each chosen pair $(S,S')$, let $S=\{s_1,\ldots,s_{k-1}\}$ and $S'=\{t_1,\ldots,t_{k-1}\}$. For each $i\in[k-1]$, choose $A_i\in N_H(s_i)\cap N_H(t_i)\cap\binom Q{k-1}$ disjoint from $S\cup S'\cup\{v,w\}\cup\bigcup_{j<i}A_j$. Note that at each step there are at most $2+2(k-1)+(i-1)(k-1)\leq k^2$ vertices to avoid. Thus there are at least $h\binom{n-1}{k-1}-k^2\binom{n-2}{k-2}\geq(h/2)\binom{n-1}{k-1}$ choices for $A_i$.
	Let $X=S\cup S'\cup\bigcup_{i=1}^{k-1}A_i$. Then $|X|=k^2-1$, and the two matchings $\{S\cup\{v\}\}\cup\{A_i\cup\{t_i\}:i\in[k-1]\}$ and $\{S'\cup\{w\}\}\cup\{A_i\cup\{s_i\}:i\in[k-1]\}$ cover $X\cup\{v\}$ and $X\cup\{w\}$, respectively.
	In total, there are at least $(\alpha/4)^2(h/2)^{k-1}\binom{n-1}{k-1}^{k+1}$ choices of $(S,S',A_1,\ldots,A_{k-1})$. Each fixed $X\subseteq V(H)$of size $k^2-1$ has at most $(k^2-1)!/((k-1)!)^{k+1}$ partitions into the sets $S,S',A_1,\ldots,A_{k-1}$. Since $\binom{n-1}{k-1}\geq n^{k-1}/(2(k-1)!)$ for sufficiently large $n$, there are at least $(\alpha/8)^2(h/4)^{k-1}n^{k^2-1}/(k^2-1)!$ distinct sets $X$. Since $\tau/4\leq\min\{\alpha/8,h/4\}<1$ and $k+1\leq2k-1$,  $(\alpha/8)^2(h/4)^{k-1}n^{k^2-1}/(k^2-1)!\geq (\tau/4)^{2k-1}n^{k^2-1}/(k^2-1)!=\beta n^{k^2-1}$.
	
	Apply Lemma~\ref{lem:LMabsorb} with $i=k$ and $\beta$. Then there exists $A\subseteq V(H)$ with $|A|\leq(\beta/2)^kn/(4k^2(k-1))$ and $|A|\in k\mathbb{N}$ such that $H[A\cup W]$ has a perfect matching for every $W\subseteq V(H)\setminus A$ satisfying $|W|\leq(\beta/2)^{2k}n/(32k^3(k-1)^2)$ and $|W|\in k\mathbb{N}$. The proof of this case is done.
\end{proof}

\section{Almost perfect matchings}\label{sec:endpoint}
We first obtain almost-perfect matchings in the non-extremal case.

\begin{definition}
Let $n,k$ be integers and $\eta$ be a constant such that $k\ge3$, $n\in k\mathbb N$, and $\eta>0$. An $n$-vertex $k$-graph $H$ is \emph{$\eta$-extremal} if there is a set $U\subseteq V(H)$ of size $n-n/k+1$ such that $e(H[U])\le\eta n^k$
and all but at most $\eta n$ vertices $u\in U$ satisfy $|d_H(u)-D(n,k,n/k)|\le\eta n^{k-1}$.
\end{definition}


The following theorem guarantees an almost perfect matching
in the non-extremal case under the stated minimum-degree
and Ore-degree conditions.
We argue by contradiction, applying the approximate
fractional cover lemma to a sequence of counterexamples
and using compactness to obtain a limiting weight distribution. The degree conditions force equality
in Theorem~\ref{thm:tail} after normalization of the limiting
weights. Its equality case yields a large sparse vertex set,
and a further application of the Ore-degree condition gives
the required degree estimates, contradicting non-extremality.
\begin{theorem}[Almost-perfect matchings in the non-extremal case]\label{thm:alternative}
	Let $k$ be an integer and $\alpha,\eta$ be constants such that $k\ge3$ and  $\alpha,\eta>0$. There exists $\delta_*=\delta_*(k,\alpha,\eta)>0$ such that, for every $0<\tau<1/2$, there is $n_0=n_0(k,\alpha,\eta,\tau)$ such that the following holds whenever $n\ge n_0$ and $n\in k\mathbb{N}$. Suppose that $H$ is an $n$-vertex $k$-graph. If $H$ is not $\eta$-extremal, $\delta_1(H)\ge\alpha\binom{n-1}{k-1}$ and 
	$\sigma_k(H)\ge\left(k-k\left(1-\frac1k\right)^{k-1}-\delta_*\right)\binom{n-1}{k-1}$,
	then $H$ has a matching covering all but at most $\tau n$ vertices.
\end{theorem}

\begin{proof}
	Since every $\eta'$-extremal hypergraph is $\eta$-extremal
whenever $0<\eta'\le\eta$, it suffices to prove the theorem
for $0<\eta\le1/(8k)$.
Suppose, on the contrary, that for some fixed $k,\alpha,\eta$ and every $\delta_*>0$, there exists $\tau\in(0,1/2)$ such that for arbitrarily large $n$ divisible by $k$, there is a non-$\eta$-extremal $n$-vertex $k$-graph $H$ satisfying the degree conditions in the theorem and $\nu(H)<\frac{n-\tau n}{k}$. Hence, for each integer $j>2$, we may choose $\tau_j\in(0,1/2)$ such that counterexamples of arbitrarily large order divisible by $k$ exist for $\delta_*=1/j$. Let $m_j$ be an integer such that $m_j\ge\max\{j,1/\tau_j\}$. Thus $m_j\to\infty$ as $j\to\infty$. We may then choose a sequence of integers $n_j\in k\mathbb{N}$ with $n_j\to\infty$ such that, for each $j$, $n_j$ is sufficiently large for the application of Lemma~\ref{lem:approx} below and there is a non-$\eta$-extremal $n_j$-vertex $k$-graph $H_j$ satisfying
	$\nu(H_j)<\frac{1-\tau_j}{k}n_j$,
	$\delta_1(H_j)\ge\alpha\binom{n_j-1}{k-1}$, and
	$\sigma_k(H_j)\ge\left(k-k\left(1-\frac1k\right)^{k-1}-\frac1j\right)\binom{n_j-1}{k-1}$.

For a weight function $z:V(H)\to[0,1]$, we write $\mathcal R(H,z)=\left\{e\in E(H):\sum_{v\in e}z(v)<1\right\}$. Applying Lemma~\ref{lem:approx} to $H_j$ with $x=\xi=(1-\tau_j)/k$ and $\rho=1/(4k^2m_j^2)$, we obtain weights $w_j:V(H_j)\to[0,1]$ such that
		$\frac1{n_j}\sum_{v\in V(H_j)}w_j(v)\le\frac{1-\tau_j}{k}+\frac1{4k^2m_j^2}$ and
	$ |\mathcal R(H_j,w_j)|\le\frac{n_j^k}{4k^2m_j^2}$.

	Let $F_j$ be the set of vertices incident with more than $n_j^{k-1}/(2km_j)$ edges of $\mathcal R(H_j,w_j)$. Since $|F_j|n_j^{k-1}/(2km_j)\leq k|\mathcal R(H_j,w_j)|$, $|F_j|\le n_j/(2m_j)$. Let $V'_j=V(H_j)\setminus F_j$. Note that $\tau_j> 1/(2m_j)+1/(4km_j^2)$. Thus 
	\begin{equation}\label{eq:meanbelow}
		\frac{\sum_{v\in V'_j}w_j(v)}{|V'_j|}\leq\frac{(1-\tau_j)/k+1/(4k^2m_j^2)}{1-1/(2m_j)}<\frac1k.
	\end{equation}
	
	\smallskip\noindent
	\textbf{Claim 1.} \emph{There is a weight function $z_j:V(H_j)\to[0,1]$ such that $\min_vz_j(v)=0$,
		$\frac1{n_j}\sum_vz_j(v)\le\frac1k+\frac1{2m_j}$ and $\max_vd_{\mathcal R(H_j,z_j)}(v)\le\frac{n_j^{k-1}}{2km_j}$.}

    \smallskip\noindent
	\emph{Proof of Claim 1.}
    Let $G_j$ be the $k$-graph on $V'_j$ with
	$E(G_j)=\left\{e\in E(H_j[V'_j]):\sum_{v\in e}w_j(v)\ge1\right\}$.
	Then $w_j|_{V'_j}$ is a fractional vertex cover of $G_j$. Let $z_j:V'_j\to[0,1]$ be a minimum fractional vertex cover of $G_j$, and let $\varphi_j:E(G_j)\to[0,1]$ be a maximum fractional matching of $G_j$. By linear programming duality and \eqref{eq:meanbelow},
	\begin{equation}\label{eq:optimalcover}
    \sum_{e\in E(G_j)}\varphi_j(e)=\sum_{v\in V'_j}z_j(v)\le\sum_{v\in V'_j}w_j(v)<\frac{|V'_j|}{k}.
    \end{equation}
	It follows that 
	\begin{equation*}
	\sum_{v\in V'_j}\sum_{\substack{e\in E(G_j)\\v\in e}}\varphi_j(e)
	=k\sum_{e\in E(G_j)}\varphi_j(e)<|V'_j|.
    \end{equation*}
	Thus there exists a vertex $v_0\in V'_j$ such that $\sum_{e\ni v_0}\varphi_j(e)<1$. We claim that $z_j(v_0)=0$. Indeed, since 
	\[
	\begin{aligned}
		0=\sum_{v\in V'_j}z_j(v)-\sum_{e\in E(G_j)}\varphi_j(e)=\sum_{v\in V'_j}z_j(v)
		\left(1-\sum_{\substack{e\in E(G_j)\\v\in e}}\varphi_j(e)\right)
		+\sum_{e\in E(G_j)}\varphi_j(e)
		\left(\sum_{v\in e}z_j(v)-1\right).
	\end{aligned}
	\]
	 and every summand on the right is nonnegative, we can derive that $z_j(v_0)(1-\sum_{\substack{v_0\in e}}\varphi_j(e))=0$. Thus $z_j(v_0)=0$.
	
	Extend $z_j$ to $V(H_j)$ by defining $z_j(v)=1$ for $v\in F_j$. By \eqref{eq:optimalcover} and $|F_j|\le n_j/(2m_j)$,
	\[
	\frac1{n_j}\sum_{v\in V(H_j)}z_j(v)
	<\frac1{n_j}(\frac{|V'_j|}{k}+|F_j|)
	=\frac1{n_j}(\frac{n_j}{k}+\left(1-\frac1k\right)|F_j|)\le\frac1k+\frac1{2m_j}.
	\]
	
	Every edge meeting $F_j$ contains a vertex $x$ with $z_j(x)=1$, and every edge of $G_j$ satisfies $\sum_{v\in e}z_j(v)\ge1$. Hence
	$\mathcal R(H_j,z_j)\subseteq\{e\in\mathcal R(H_j,w_j):e\subseteq V'_j\}$.
	For $v\in V'_j$, the definition of $F_j$ therefore gives
	$d_{\mathcal R(H_j,z_j)}(v)\le d_{\mathcal R(H_j,w_j)}(v)\le n_j^{k-1}/(2km_j)$,
	whereas $d_{\mathcal R(H_j,z_j)}(v)=0$ for $v\in F_j$.
	This completes the proof of Claim 1.
    \claimqed
    \smallskip


	Let $\mu_j$ be the discrete probability measure on $[0,1]$ defined by
	$\mu_j(\{t\})=\frac{|\{v\in V(H_j):z_j(v)=t\}|}{n_j}$ for each $t\in[0,1]$.
	Thus $\mu_j$ is the empirical probability measure of the weights.
	The space of probability measures on $[0,1]$ is weakly compact. Passing to a subsequence and retaining the index $j$, we may assume that $\mu_j$ converges
	weakly to a probability measure $\mu$ on $[0,1]$. 
	
	\medspace
	Since $\mu_j$ converges weakly to  $\mu$,  $\mu_j^{\otimes(k-1)}$ converges
	weakly to $\mu^{\otimes(k-1)}$, where $\mu_j^{\otimes(k-1)}$ is the joint distribution of $k-1$
	independent and identically distributed random variables with distribution $\mu_j$. Let $Z$ be a random variable with distribution $\mu$. Since the identity function $t\mapsto t$ is bounded and continuous on $[0,1]$, 
	\[
	\E Z
	=\int_{[0,1]}t\,d\mu(t)
	=\lim_{j\to\infty}\int_{[0,1]}t\,d\mu_j(t)
	=\lim_{j\to\infty}\frac1{n_j}
	\sum_{v\in V(H_j)}z_j(v).
	\]
	Thus $\E Z\leq 1/k$ by Claim 1.
	Let $b=\inf\supp(\mu)$, where $$\operatorname{supp}(\mu)
=
\left\{
x\in[0,1]:
\mu\bigl((x-\varepsilon,x+\varepsilon)\cap[0,1]\bigr)>0
\text{ for every }\varepsilon>0
\right\}.$$ Then $\Prob(Z\ge b)=1$ and $\Prob(Z<b+1/\ell)>0$ for every positive integer $\ell$. Together with the mean bound above, this implies $0\le b\le\E Z\le1/k$.
	
	\smallskip\noindent
	\textbf{Claim 2.} \emph{There are sets $T_j\subseteq V(H_j)$ of size $n_j/k-1$ such that, with $U_j=V(H_j)\setminus T_j$, we have $e(H_j[U_j])\le \eta^{2k} n_j^k$ for all sufficiently large $j$.}
	
	\smallskip\noindent
	\emph{Proof of Claim 2.}
	For each $j$ and arbitrary fixed $x_j\in V(H_j)$, choose vertices $U_{j,1},\ldots,U_{j,k-1}$ independently and uniformly from $V(H_j)$, with replacement. Thus $z_j(U_{j,1}),\ldots,z_j(U_{j,k-1})$ are independent and identically distributed random variables with distribution $\mu_j$.
	
	There are exactly $d_{H_j}(x_j)-d_{\mathcal R(H_j,z_j)}(x_j)$ edges $e\in E(H_j)$ containing $x_j$ and satisfying $\sum_{v\in e}z_j(v)\ge1$. Hence
	\[
	\begin{aligned}
		\Prob\bigl(\{x_j,U_{j,1},\ldots,U_{j,k-1}\}
		\in E(H_j)\setminus\mathcal R(H_j,z_j)\bigr)
		=\frac{(k-1)!}{n_j^{k-1}}
		\left(d_{H_j}(x_j)-d_{\mathcal R(H_j,z_j)}(x_j)\right).
	\end{aligned}
	\]
	By the definition of $\mathcal R(H_j,z_j)$, the event on the left implies $z_j(x_j)+\sum_{i=1}^{k-1}z_j(U_{j,i})\ge1$. Therefore,
	\begin{equation}\label{eq:coveredlinktail}
		\frac{(k-1)!}{n_j^{k-1}}
		\left(d_{H_j}(x_j)-d_{\mathcal R(H_j,z_j)}(x_j)\right)
		\le\Prob\left(z_j(x_j)+\sum_{i=1}^{k-1}z_j(U_{j,i})\ge1\right).
	\end{equation}
	By Claim 1, $d_{\mathcal R(H_j,z_j)}(x_j)\le n_j^{k-1}/(2km_j)$ for every $x_j\in V(H_j)$. Combining with inequality \eqref{eq:coveredlinktail},
	$$\frac{(k-1)!}{n_j^{k-1}}
	d_{H_j}(x_j)
	\le\Prob\left(z_j(x_j)+\sum_{i=1}^{k-1}z_j(U_{j,i})\ge1\right)+\frac{(k-1)!}{2km_j}.$$
	
	Since $d_{H_j}(x_j)/\binom{n_j-1}{k-1}\le1$ and $0\le1-\frac{(k-1)!\binom{n_j-1}{k-1}}{n_j^{k-1}}=1-\prod_{i=1}^{k-1}\left(1-\frac{i}{n_j}\right)\le\frac{k(k-1)}{2n_j}$, one can see that
	\begin{equation}\label{eq:finitevertexprob}
		\begin{split}
			\frac{d_{H_j}(x_j)}{\binom{n_j-1}{k-1}}&=\frac{(k-1)!d_{H_j}(x_j)}{n_j^{k-1}}+\left(1-\frac{(k-1)!\binom{n_j-1}{k-1}}{n_j^{k-1}}\right)\frac{d_{H_j}(x_j)}{\binom{n_j-1}{k-1}}\\
			&\le\Prob\left(z_j(x_j)+\sum_{i=1}^{k-1}z_j(U_{j,i})\ge1\right)
			+\frac{(k-1)!}{2km_j}+\frac{k(k-1)}{2n_j}.
		\end{split}
	\end{equation}
	
	Let $a\in[0,1]$ be fixed, and suppose that $x_j\in V(H_j)$ satisfies $z_j(x_j)\le a$ for all sufficiently large $j$. For every sufficiently large $j$, $z_j(x_j)+\sum_{i=1}^{k-1}z_j(U_{j,i})\ge1$
	implies $a+\sum_{i=1}^{k-1}z_j(U_{j,i})\ge1$.
	Hence by \eqref{eq:finitevertexprob},
	\[
	\frac{d_{H_j}(x_j)}{\binom{n_j-1}{k-1}}
	\le
	\Prob\left(a+\sum_{i=1}^{k-1}z_j(U_{j,i})\ge1\right)
	+\frac{(k-1)!}{2km_j}+\frac{k(k-1)}{2n_j}.
	\]
	The random vector $(z_j(U_{j,1}),\ldots,z_j(U_{j,k-1}))$ has distribution $\mu_j^{\otimes(k-1)}$. Since $a$ is fixed, the set $\left\{(t_1,\ldots,t_{k-1})\in[0,1]^{k-1}:a+t_1+\cdots+t_{k-1}\ge1\right\}$
	is a fixed closed subset of $[0,1]^{k-1}$. Therefore, by the fact that
	$\mu_j^{\otimes(k-1)}$ converges weakly to $\mu^{\otimes(k-1)}$ and Portmanteau theorem~\cite[Theorem~2.1(iii)]{Billingsley},
	\begin{equation}\label{eq:sampling}
		\limsup_{j\to\infty}
		\frac{d_{H_j}(x_j)}{\binom{n_j-1}{k-1}}
		\le
		\Prob(a+Z_1+\cdots+Z_{k-1}\ge1),
	\end{equation}
	where $Z_1,\ldots,Z_{k-1}$ are independent copies of $Z$.
	
	By Claim 1, there exists a vertex $v_j\in V(H_j)$ such that $z_j(v_j)=0$. Suppose that $b=1/k$. Since $\Prob(Z\ge b)=1$ and $\E Z\le1/k$, we have $\Prob(Z=1/k)=1$. Thus $\Prob(Z_1+\cdots+Z_{k-1}\ge1)=0$. Applying \eqref{eq:sampling} to $v_j$ with $a=0$, 
	\[
	\begin{aligned}
		\liminf_{j\to\infty}
		\frac{d_{H_j}(v_j)}{\binom{n_j-1}{k-1}}
		\le\limsup_{j\to\infty}
		\frac{d_{H_j}(v_j)}{\binom{n_j-1}{k-1}}\le
		\Prob(Z_1+\cdots+Z_{k-1}\ge1)=0,
	\end{aligned}
	\]
	contradicting the condition that $\delta_1(H_j)\ge\alpha\binom{n_j-1}{k-1}$. Since $b\le1/k$, it follows that $b<1/k$.

	Let  $\ell$ be a sufficiently large positive integer such  that $k(b+1/\ell)<1$, and let
	$B_j=\{v\in V(H_j):b-1/\ell<z_j(v)<b+1/\ell\}$.
	The interval $(b-1/\ell,b+1/\ell)$ is an open neighborhood of $b\in\supp(\mu)$, so $\mu(b-1/\ell,b+1/\ell)>0$. Let $\mu_0=\mu(b-1/\ell,b+1/\ell)$. By weak convergence and the Portmanteau theorem,
	\[
	\begin{aligned}
		\liminf_{j\to\infty}\frac{|B_j|}{n_j}=\liminf_{j\to\infty}\mu_j((b-1/\ell,b+1/\ell))
		\ge\mu_0>0.
	\end{aligned}
	\]
	Consequently, $|B_j|\ge \mu_0 n_j/2$ for all sufficiently large $j$. For every $e\in\binom{B_j}{k}$, we have $\sum_{v\in e}z_j(v)<k(b+1/\ell)<1$. Thus every edge of $H_j[B_j]$ belongs to $\mathcal R(H_j,z_j)$. By Claim 1,
	\begin{equation}\label{eq:exceptiontotal}
		|\mathcal R(H_j,z_j)|
		=\frac1k\sum_{v\in V(H_j)}d_{\mathcal R(H_j,z_j)}(v)
		\le\frac{n_j^k}{2k^2m_j}.
	\end{equation}
	For fixed $\ell$ and all sufficiently large $j$,
	$\binom{|B_j|}{k}
	\ge\frac{(\mu_0 n_j/4)^k}{k!}
	>\frac{n_j^k}{2k^2m_j}
	\ge|\mathcal R(H_j,z_j)|$,
	Hence there exists a nonedge $e_j=\{y_{j,1},\ldots,y_{j,k}\}\in\binom{B_j}{k}\setminus E(H_j)$.
	
	By the degree condition,
	$\frac{\sum_{i=1}^k d_{H_j}(y_{j,i})}
	{k\binom{n_j-1}{k-1}}
	\ge1-\left(1-\frac1k\right)^{k-1}-\frac1{kj}$.
	Since $z_j(y_{j,i})<b+1/\ell$ and $\ell$ is fixed, we may apply \eqref{eq:sampling} to each sequence $y_{j,i}$ with $a=b+1/\ell$. It follows that
	\[
	\begin{aligned}
		1-\left(1-\frac1k\right)^{k-1}
		\le\frac1k\sum_{i=1}^k\limsup_{j\to\infty}
		\frac{d_{H_j}(y_{j,i})}{\binom{n_j-1}{k-1}}
		\le\Prob\left(b+\frac1\ell+Z_1+\cdots+Z_{k-1}\ge1\right)
	\end{aligned}
	\]
 holds for every sufficiently large integer $\ell$. By continuity of probability, we obtain
	\begin{equation}\label{eq:endpointanchor}
		1-\left(1-\frac1k\right)^{k-1}
		\le\Prob(b+Z_1+\cdots+Z_{k-1}\ge1).
	\end{equation}
	
	Recall that $b<1/k$. Let $Y_i=(Z_i-b)/(1/k-b)$ for each $i\in[k-1]$. Thus $Y_1,\ldots,Y_{k-1}$ are independent and nonnegative random variables with
	$\E Y_i=\frac{\E Z-b}{1/k-b}\le1$.
	Note that $b+\sum_{i=1}^{k-1}Z_i=kb+(1/k-b)\sum_{i=1}^{k-1}Y_i$. Applying Theorem~\ref{thm:tail} with $r=k-1$ and $T=k$, we can derive that
	\[
	\begin{aligned}
		\Prob(b+Z_1+\cdots+Z_{k-1}\ge1)
		=\Prob(Y_1+\cdots+Y_{k-1}\ge k)
		\le1-\left(1-\frac1k\right)^{k-1}.
	\end{aligned}
	\]
	Together with \eqref{eq:endpointanchor}, this forces equality. The equality statement in Theorem~\ref{thm:tail} implies that each $Y_i$ takes values $0$ and $k$ with probabilities $1-1/k$ and $1/k$, respectively. Consequently $Z$ takes values $b$ and $1-(k-1)b$ with probabilities $1-1/k$ and $1/k$, respectively. Note that $1-(k-1)b-1/k=(k-1)(1/k-b)>0$. 
	It follows that $b<1/k<1-(k-1)b$.
	
	We regard $\mu_j$ and $\mu$ as probability measures on $\mathbb R$
	by assigning zero mass outside $[0,1]$. This identification preserves
	weak convergence. Since $b<1/k<1-(k-1)b$, we have
	$\mu((-\infty,1/k))=\mu((-\infty,1/k])=1-\frac1k$.
	Let $U_j^0=\{v\in V(H_j):z_j(v)<1/k\}$. By the definition of
	$\mu_j$, we have $|U_j^0|/n_j=\mu_j((-\infty,1/k))$.
	Since $(-\infty,1/k)$ is open and $(-\infty,1/k]$ is closed,
	weak convergence and the Portmanteau theorem~\cite[Theorem~2.1(iv)]{Billingsley} imply that
	\[
	\begin{aligned}
		1-\frac1k&=\mu((-\infty,1/k))
		\le\liminf_{j\to\infty}\mu_j((-\infty,1/k))
		=\liminf_{j\to\infty}\frac{|U_j^0|}{n_j}\\
		&\le\limsup_{j\to\infty}\frac{|U_j^0|}{n_j}
		=\limsup_{j\to\infty}\mu_j((-\infty,1/k))
		\le\limsup_{j\to\infty}\mu_j((-\infty,1/k])
		\le\mu((-\infty,1/k])=1-\frac1k.
	\end{aligned}
	\]
	It follows that $(1-1/k)n_j-\eta^{2k} n_j/2\leq |U_j^0|\leq(1-1/k)n_j+\eta^{2k} n_j/2$ for sufficiently large $j$.
	
	Since every edge $e$ of $H_j[U_j^0]$ satisfies $\sum_{v\in e}z_j(v)<k(1/k)=1$, 
	$E(H_j[U_j^0])\subseteq\mathcal R(H_j,z_j)$ and
	$e(H_j[U_j^0])\le|\mathcal R(H_j,z_j)|\le\frac{n_j^k}{2k^2m_j}$ by \eqref{eq:exceptiontotal}.
	
	Let $U_j$  be a subset of $V(H_j)$  such that $|U_j|=n_j-n_j/k+1$ and
	$|U_j\cap U_j^0|\geq (1-1/k)n_j-\eta^{2k} n_j/2$.
	Therefore,
	\[
	\begin{aligned}
	e(H_j[U_j])\le e(H_j[U_j^0])+|U_j\setminus U_j^0|\binom{n_j-1}{k-1}
		\le \frac{n_j^k}{2k^2m_j}+\binom{n_j-1}{k-1}(\eta^{2k} n_j/2+1)\leq \eta^{2k} n_j^k
	\end{aligned}
	\]
	for sufficiently large $j$. This proves Claim 2.
    \claimqed
    \smallskip

	\smallskip\noindent
	\textbf{Claim 3.} \emph{For all sufficiently large $j$, all but at most $\eta n_j$ vertices $v\in U_j$ satisfy $|d_{H_j}(v)-D(n_j,k,n_j/k)|\le\eta n_j^{k-1}$.}
	
	\smallskip\noindent
	\emph{Proof of Claim 3.}
	Let $D_j=D(n_j,k,n_j/k)$. Since
	$\frac{\binom{n_j-n_j/k}{k-1}}{\binom{n_j-1}{k-1}}
	=\prod_{i=0}^{k-2}
	\frac{n_j-n_j/k-i}{n_j-1-i}
	\ge\left(1-\frac1k\right)^{k-1}$,
	we have
	$D_j\le
	\left(1-\left(1-\frac1k\right)^{k-1}\right)
	\binom{n_j-1}{k-1}$.
	The degree-sum condition on $H_j$ implies
	\[
	\begin{aligned}
		\sigma_k(H_j)\ge\left(k-k\left(1-\frac1k\right)^{k-1}-\frac1j\right)
		\binom{n_j-1}{k-1}\ge kD_j-\frac1j\binom{n_j-1}{k-1}.
	\end{aligned}
	\]
	
	By Claim 2, for all sufficiently large $j$ we have
	$e(H_j[U_j])<\eta^{2k} n_j^k$.
	 Let $U_j'=\{v\in U_j:|d_{H_j}(v)-D_j|> \eta n_j^{k-1}\}$, $U_j^1=\{v\in U_j:d_{H_j}(v)<D_j-\eta n_j^{k-1}\}$ and $U_j^2=\{v\in U_j:d_{H_j}(v)>D_j+\eta n_j^{k-1}\}$. Thus $U_j^1\dot\cup U_j^2=U_j'$. It suffices to prove that $|U_j^1|+|U_j^2|\leq \eta n_j$.
	
	For arbitrary $k$ distinct vertices $\{x_1,\ldots,x_k\}\subseteq U_j^1$, $\sum_{i\in [k]}d_{H_j}(x_i)< kD_j-k\eta n_j^{k-1}\le\sigma_k(H_j)$ for sufficiently large $j$. 
	Thus every $k$-set in $U_j^1$ is an edge of $H_j$.
	Suppose that $|U_j^1|\ge\eta n_j/2$. Then
	\[
	\begin{aligned}
		e(H_j[U_j])\ge\binom{|U_j^1|}{k}\ge\frac{(|U_j^1|-k+1)^k}{k!}\ge\frac{(\eta n_j/4)^k}{k!}
		=\frac{\eta^k}{4^k k!}n_j^k> \eta^{2k} n_j^{k},
	\end{aligned}
	\]
	contradicting Claim 2. Hence $|U_j^1|<\eta n_j/2$.
	 
	For each $v\in U_j$, the number of possible edges containing $v$
	and intersecting $T_j$ is at most
	\[
	\binom{n_j-1}{k-1}-\binom{|U_j|-1}{k-1}
	=\binom{n_j-1}{k-1}-\binom{n_j-n_j/k}{k-1}
	=D_j.
	\]
	Consequently, every $v\in U_j^2$ satisfies
	$d_{H_j[U_j]}(v)
	\ge d_{H_j}(v)-D_j
	\geq\eta n_j^{k-1}$. It follows that
	\[
	|U_j^2|\eta n_j^{k-1}\leq \sum_{v\in U_j^2}d_{H_j[U_j]}(v)\leq k e(H_j[U_j])
	\leq k\eta^{2k}n_j^k.
	\]
	Hence $|U_j^2|<\eta n_j/2$ and $|U_j^1|+|U_j^2|<\eta n_j$.
	This proves Claim 3.
    \claimqed
    \smallskip
	
	For all sufficiently large $j$,   $H_j$ is $\eta$-extremal by Claim 2 and Claim 3, a contradiction.
	\end{proof}

\section{Non-extremal case}\label{sec:non-extremal}

We now combine Theorem~\ref{thm:alternative} with
Lemma~\ref{lem:absorb} to obtain a perfect matching
in the non-extremal case under the minimum-degree and the Ore-degree conditions.

\begin{lemma}\label{prop:nonextremal}
	Let $k$ be an integer and $\alpha,\eta$ be constants such that $k\ge3$ and  $\alpha,\eta>0$. There exist $\delta=\delta(k,\alpha,\eta)>0$ and $n_0=n_0(k,\alpha,\eta)$ such that the following holds for every $n\ge n_0$ with $n\in k\mathbb N$. If $H$ is an $n$-vertex $k$-graph with $\delta_1(H)\ge\alpha\binom{n-1}{k-1}$ and $\sigma_k(H)\ge kD(n,k,n/k)-\delta\binom{n-1}{k-1}$, and $H$ is not $\eta$-extremal, then $H$ has a perfect matching.
\end{lemma}

\begin{proof}
Applying Theorem~\ref{thm:alternative} with
$\alpha/2$ and $\eta/8$, and decreasing the resulting
constant if necessary, we obtain $0<\delta_*\le1$.
	Choose $\zeta$ sufficiently small and $\beta=\frac{(\alpha\zeta/(512k^2))^{2k-1}}{(k^2-1)!}$ such that
	$\zeta\le\frac{k}{2}\left(\frac12-\left(1-\frac1k\right)^{k-1}\right)$
	and 
	$\frac{(\beta/2)^k}{8k(k-1)}
	\le\min\left\{\frac{\alpha}{4k},\frac{\delta_*}{8k^2},\frac{\eta}{4}\right\}$.
	Let $\delta$ be a constant such that $0<\delta\le\min\{\delta_*/8,\zeta/2\}$.
	Let $n_0$ be a sufficiently large integer such that the conclusion of Lemma~\ref{lem:absorb} and Theorem~\ref{thm:alternative} hold for $n>n_0$ and $\frac{k(k-1)}{2(n_0-k+1)}\le\delta$.

	By the choices of $\zeta$ and $\delta$, $k-k(1-1/k)^{k-1}-2\delta\ge k/2+\zeta$. It follows that $\sigma_k(H)\geq (k/2+\zeta)\binom{n-1}{k-1}$. By Lemma~\ref{lem:absorb}, there is a  set $A$ such that
	$|A|\le\frac{(\beta/2)^k}{8k(k-1)}n$, $k\mid |A|$ and $H[A\cup W]$ has a perfect matching for every $W\subseteq V(H)\setminus A$ satisfying $|W|\le(\beta/2)^{2k}n/(32k^3(k-1)^2)$ with $k\mid |W|$.
	
	Let $G=H-A$ and $N=n-|A|$. Then $N\in k\mathbb N$ and $N\ge n/2$.
	Let $\Delta=\binom{n-1}{k-1}-\binom{N-1}{k-1}$. Note that $\Delta\le |A|\binom{n-2}{k-2}$.
	For every $v\in V(G)$, we have $0\le d_H(v)-d_G(v)\le\Delta$. Since $\binom{n-2}{k-2}=\frac{k-1}{n-1}\binom{n-1}{k-1}$, the choice of $\beta$ implies
	$\Delta\le\min\left\{\frac{\alpha}{2},\frac{\delta_*}{4k}\right\}\binom{n-1}{k-1}$,
	$|A|\le\frac{\eta n}{4}$ and $\Delta\le\frac{\eta n^{k-1}}4$.
	Consequently,
	\[
	\delta_1(G)\ge\alpha\binom{n-1}{k-1}-\Delta
	\ge\frac{\alpha}{2}\binom{n-1}{k-1}
	\ge\frac{\alpha}{2}\binom{N-1}{k-1}.
	\]
	Every non-edge of $G$ is also a non-edge of $H$. Thus by the choice of $\delta$,
	\begin{equation*}\label{eq:nonextremalafterdeletion}
		\begin{aligned}
			\sigma_k(G)
			&\ge kD(n,k,n/k)-\delta\binom{n-1}{k-1}-k\Delta\\
			&=\left(k-k\prod_{i=0}^{k-2}
			\left(1-\frac1k+\frac{k-1-i}{k(n-1-i)}\right)-\delta\right)\binom{n-1}{k-1}-k\Delta\\
			&\ge\left(k-k\left(1-\frac1k\right)^{k-1}-\delta_*\right)\binom{N-1}{k-1}.
		\end{aligned}
	\end{equation*}

	Next we show that $G$ is not $(\eta/8)$-extremal. Otherwise, there exists $U_G\subseteq V(G)$ with $|U_G|=N-N/k+1$ such that $e(G[U_G])\le(\eta/8)N^k$ and all but at most $\eta N/8$ vertices $u\in U_G$ satisfy $|d_G(u)-D(N,k,N/k)|\le\frac{\eta}{8}N^{k-1}$.
	
	Choose $A_T\subseteq A$ with $|A_T|=|A|-|A|/k$, and let $U_H=U_G\cup A_T$. Every edge of $H[U_H]$ that is not an edge of $G[U_G]$ intersects $A$, so
	\begin{equation}\label{eq:transfersparse}
		e(H[U_H])\le\frac{\eta}{8}N^k+|A|\binom{n-1}{k-1}
		\le\frac{3\eta}{8}n^k<\eta n^k.
	\end{equation}
	
	Note that $0\le D(n,k,n/k)-D(N,k,N/k)\le\Delta$. For $u\in U_G$, $0\le d_H(u)-d_G(u)\le \Delta$. Thus
	\[
	|d_H(u)-D(n,k,n/k)|
\le |d_G(u)-D(N,k,N/k)|+\Delta
	\]
	for each $u\in U_G$.
	Hence all but at most $\eta N/8$ vertices $u\in U_G$ satisfy
	\[
	|d_H(u)-D(n,k,n/k)|
	\le\frac{\eta}{8}N^{k-1}+\Delta
	\le\frac{3\eta}{8}n^{k-1}<\eta n^{k-1}.
	\]
	The remaining vertices of $U_H$ number at most $\eta N/8+|A_T|\le3\eta n/8<\eta n$. Together with \eqref{eq:transfersparse}, this shows that $H$ is $\eta$-extremal, a contradiction.
	
	Apply Theorem~\ref{thm:alternative} to $G$ with $k,\alpha/2,\eta/8$ and $\tau=(\beta/2)^{2k}/(32k^3(k-1)^2)$. We obtain a matching $M'$ in $G$ whose uncovered set $W=V(G)\setminus V(M')$ satisfies
	\[
	|W|\le\frac{(\beta/2)^{2k}}{32k^3(k-1)^2}N
	\le\frac{(\beta/2)^{2k}}{32k^3(k-1)^2}n.
	\]
	Since $N\in k\mathbb N$, we have $k\mid |W|$. By the absorption property of $A$, there is a perfect matching of $H[A\cup W]$. Together with $M'$, this forms a perfect matching of $H$.
\end{proof}

\section{Extremal case}\label{sec:extremal}

In this section, we handle the extremal case. For convenience
in later applications, we prove the result for
$n=ks-1+c$, where $c\ge1$ is fixed, including the
perfect matching case $c=1$.
We assume that $H-Z_0$ is $\eta$-extremal for some
$Z_0\subseteq V(H)$ of size $c-1$, while imposing the strict
Ore-degree condition on $H$.
Using the extremal structure of $H-Z_0$, we construct an
$s$-matching in $H$.

\begin{lemma}\label{lem:localcover}
	Let $k,c$ be integers such that $k\ge3$ and $c\ge1$. There exist $\eta>0$ and $s_0$ such that the following holds for integers $n,s$ with $s\ge s_0$ and $n=ks-1+c$. Suppose that $H$ is an $n$-vertex $k$-graph satisfying $\sigma_k(H)>kD(n,k,s)$. If $H-Z_0$ is $\eta$-extremal for some $Z_0\subseteq V(H)$ of size $c-1$, then $\nu(H)\ge s$.
\end{lemma}
\begin{proof}
	Fix $0<\eta<(2k)^{-2k-4}$, and choose
$s_0=s_0(k,c,\eta)$ sufficiently large.
	Let $G=H-Z_0$ and write $N=n-c+1=ks$. Since $G$ is $\eta$-extremal, there exists $U_0\subseteq V(G)$ of size $N-N/k+1$ such that $e(H[U_0])\le\eta N^k$
	and all but at most $\eta N$ vertices $u\in U_0$ satisfy $|d_G(u)-D(N,k,N/k)|\le\eta N^{k-1}$. Let $T=V(G)\setminus U_0$. Denote the set of vertices $u\in U_0$ satisfying $|d_G(u)-D(N,k,N/k)|\le\eta N^{k-1}$ by $U'$. Thus $|U_0\setminus U'|\leq \eta N$ and for every $u\in U'$, $$|N_{H_k(U_0,T)}(u)\setminus N_G(u)|=D(N,k,N/k)-(d_G(u)-d_{G[U_0]}(u))\leq \eta N^{k-1}+d_{G[U_0]}(u).$$
	 Since every edge of $H_k(U_0,T)$ intersects $U_0$, we obtain
	\begin{equation}\label{eq:recovermixed}
		\begin{aligned}
			|E(H_k(U_0,T))\setminus E(G)|
			&\le\sum_{u\in U_0}|N_{H_k(U_0,T)}(u)\setminus N_G(u)|\\
			&=\sum_{u\in U'}|N_{H_k(U_0,T)}(u)\setminus N_G(u)|+\sum_{u\in U_0\setminus U'}|N_{H_k(U_0,T)}(u)\setminus N_G(u)|\\
			&\le\eta N^k+ke(G[U_0])+\eta N D(N,k,N/k)\le(k+2)\eta N^k.
		\end{aligned}
	\end{equation}
	Define $X=Z_0\cup\bigl\{v\in V(G):
		|N_{H_k(U_0,T)}(v)\setminus N_G(v)|
		>(2k)^{k+4}\eta N^{k-1}\bigr\}
		\cup\bigl\{u\in U_0:
		d_G(u)>D(N,k,N/k)+\eta N^{k-1}\bigr\}$
	and $b=|X|$. Hence, by \eqref{eq:recovermixed} and $|Z_0|<\eta n$,
	\begin{equation}\label{eq:badbound}
		\begin{aligned}
			b\le |Z_0|+\eta N+
			\frac{k|E(H_k(U_0,T))\setminus E(G)|}
			{(2k)^{k+4}\eta N^{k-1}}\le\left(\frac{k(k+2)}{(2k)^{k+4}}+2\eta\right)n<\frac{(k+1)^2}{(2k)^{k+4}}n
			\le\frac{n}{2^{k+4}k^2}.
		\end{aligned}
	\end{equation}
	 
Let $Q=U_0\setminus X$.	
For every $u\in Q$, 
	\begin{equation}\label{eq:Qupper}
		\begin{aligned}
			d_H(u)
			&\le d_G(u)+(c-1)\binom{n-2}{k-2}\\
			&\le D(N,k,N/k)+\eta N^{k-1}+(c-1)n^{k-2}\le D(n,k,s)+2\eta n^{k-1}.
		\end{aligned}
	\end{equation}
	Let $S\subseteq T\setminus X$ be a set of size $|S|=s-1-b$. Let $J=H-S$ and $N'=|V(J)|=n-s+1+b$.  Thus
	\begin{equation}\label{eq:Jore}
		\begin{aligned}
			\sigma_k(J)
			&>kD(n,k,s)-k\left(\binom{n-1}{k-1}-\binom{N'-1}{k-1}\right)\\
			&=k\left(\binom{N'-1}{k-1}-\binom{n-s}{k-1}\right)=kD(N',k,b+1),
		\end{aligned}
	\end{equation}
	since $N'-(b+1)=n-s$. If $b\ge1$, then
	$N'\ge n-s+1\ge\frac{k-1}{k}n>3k^2b$.
	By Theorem~\ref{thm:bpr}, there is a matching of size $b+1$ in $J$. If $b=0$, then there is an edge in $J$ by \eqref{eq:Jore}. Among all $(b+1)$-matchings of $J$, we choose $M_0$ maximizing $|V(M_0)\cap X|$.
	
	\smallskip\noindent
		\textbf{Claim 4.} \emph{There is a matching $M_1$ in $H-V(M_0)$ such that $X\subseteq V(M_0\cup M_1)$ and every edge of $M_1$ contains exactly one vertex of $S$.}
        
	\smallskip\noindent
\emph{Proof of Claim 4.}
If $X\subseteq V(M_0)$, take $M_1=\varnothing$. Otherwise, let
$v\in X\setminus V(M_0)$. Since $|M_0|=b+1>|X|$, some edge
$f\in M_0$ avoids $X$. By \eqref{eq:badbound} and $n-s=((k-1)n+c-1)/k$,
$|Q\setminus V(M_0)|\ge N-s+1-b-k(b+1)\ge k-1$. For each $Y_v\in\binom{Q\setminus V(M_0)}{k-1}$, we claim that $\{v\}\cup Y_v\notin E(J)$. Indeed, if $\{v\}\cup Y_v\in E(J)$, then $M_0'=(M_0\setminus \{f\})\cup \{\{v\}\cup Y_v\}$ is a matching of size $b+1$ satisfying
$|V(M_0')\cap X|>|V(M_0)\cap X|$, contrary to the choice of $M_0$.
Thus $\{v\}\cup Y_v\notin E(H)$ and by
\eqref{eq:Qupper},
\[
kD(n,k,s)<\sigma_k(H)
\le d_H(v)+\sum_{u\in Y_v}d_H(u)
\le d_H(v)+(k-1)D(n,k,s)+2(k-1)\eta n^{k-1}.
\]
Since $v$ is arbitrary, we conclude that
\begin{equation}\label{eq:localuncoveredlower}
	d_H(v)>D(n,k,s)-2(k-1)\eta n^{k-1}
\end{equation}
holds for every $v\in X\setminus V(M_0)$. We shall also use the following consequence of the choice of \(M_0\). For every \(v\in X\setminus V(M_0)\), every edge of \(H\) containing \(v\) and disjoint from \(V(M_0)\) must intersect \(S\). Indeed, otherwise such an edge belongs to \(J\), and replacing an edge of \(M_0\) avoiding \(X\) would increase the number of vertices of \(X\) covered by the matching, a contradiction.

We now construct $M_1$ greedily. Let $t_0=|X\setminus V(M_0)|$ and $X\setminus V(M_0)=\{v_1,v_2,\ldots,v_{t_0}\}$. Suppose that we have found
pairwise disjoint edges $e_1,\ldots,e_t$ in $H-V(M_0)$ such that
$e_i=\{v_i,x_i\}\cup Y_i$, $x_i\in S\setminus V(M_0)$ and $Y_i\in\binom {Q\setminus V(M_0)}{k-2}$ for each $i\in[t]$, where
$0\le t<t_0$.
Note that $|e_i\cap X|=1$ for each $i\in [t]$. Let
$v_{t+1}\in X\setminus\left(V(M_0)\cup\bigcup_{i=1}^t e_i\right)$.
We show that there is an edge $e_{t+1}$ containing $v_{t+1}$ such that $|e_{t+1}\cap S|=1$, $|e_{t+1}\cap Q|=k-2$ and $|e_{t+1}\cap (V(M_0)\cup\bigcup_{i=1}^t e_i)|=0$.

Let $E_{v_{t+1}}=\{e\in E(H):v_{t+1}\in e, |e\cap S|\geq 2\}$. As $|S|\leq s-1$,
\[
|E_{v_{t+1}}|\leq \binom{n-1}{k-1}-\binom{n-s}{k-1}
-(s-1)\binom{n-s}{k-2}
=D(n,k,s)-(s-1)\binom{n-s}{k-2}.
\]
Let $F_t:=\bigl(V(H)\setminus(S\cup Q\cup\{v_{t+1}\})\bigr)
\cup V(M_0)\cup\bigcup_{i=1}^{t}e_i$. Since $|V(H)\setminus(S\cup Q)|=|T\setminus S|+|X\setminus T|\le2b$ and $|V(M_0)\cup\bigcup_{i=1}^t e_i|=k(b+1+t)$, the number of edges containing $v_{t+1}$ and intersecting $F_t$ is at most $(2b+k(b+1+t))\binom{n-2}{k-2}$. Recall that every edge containing $v_{t+1}$ and avoiding $F_t$ intersects $S$. 
Therefore, there are at least
$d_H(v_{t+1})-(2b+k(b+1+t))\binom{n-2}{k-2}-|E_{v_{t+1}}|$ edges $e$ containing $v_{t+1}$ such that $|e\cap S|=1$, $|e\cap Q|=k-2$ and $|e\cap (V(M_0)\cup\bigcup_{i=1}^t e_i)|=0$. 
By \eqref{eq:localuncoveredlower} and \eqref{eq:badbound}, 
\begin{equation*}\label{eq:localgreedycount}
	\begin{aligned}
		&d_H(v_{t+1})-(2b+k(b+1+t))\binom{n-2}{k-2}-|E_{v_{t+1}}|\\
		\geq &(s-1)\binom{n-s}{k-2}-2(k-1)\eta n^{k-1}
		-\bigl(k(b+1+t)+2b\bigr)\binom{n-2}{k-2}>0
	\end{aligned}
\end{equation*}
for sufficiently large $n$.
Hence there is an edge $e_{t+1}$ containing $v_{t+1}$ such that $|e_{t+1}\cap S|=1$, $|e_{t+1}\cap Q|=k-2$ and $|e_{t+1}\cap (V(M_0)\cup\bigcup_{i=1}^t e_i)|=0$. 

Continuing this process, we obtain the desired matching \(M_1\).
This proves Claim 4.
\claimqed
\smallskip

	Let $S_0=S\setminus V(M_1)$ and $W=V(H)\setminus\bigl(S\cup V(M_0\cup M_1)\bigr)$. The uncovered vertices are precisely $S_0\mathbin{\dot\cup}W$. Note that
		$a:=|S_0|=s-1-b-|M_1|\ge s-1-2b
		\ge\frac n{2k}$ and $|W|-(k-1)a=n-k(b+1+|M_1|)-k(s-1-b-|M_1|)=c-1$.
        Note that $|W\cap T|
\le |T\setminus S|
=b<a.$
Let $W^*\subseteq W$ be a set of size $|W^*|=(k-1)a$
and $W_1\mathbin{\dot\cup}\cdots\mathbin{\dot\cup}W_{k-1}$ be a partition of $W^*$ such that
	$|W_i|=a$ for $i\in[k-1]$ and $W^*\cap T\subseteq W_1$.
Thus $W_2,\ldots,W_{k-1}\subseteq U_0$.

	Let $J$ be the complete $k$-partite $k$-graph with parts $S_0,W_1,\ldots,W_{k-1}$ and $J_0$  be the subgraph of $J$ with edge set $E(J_0)=E(J)\cap E(G)$. Since $k\ge3$, every edge of $J$ intersects both $T$ and $U_0$. Thus $E(J)\subseteq E(H_k(U_0,T))$. Moreover, every vertex of $J_0$ belongs to $V(G)\setminus X$. Hence, for every $v\in V(J_0)$,
	\[
	\begin{aligned}
		d_{J_0}(v)\ge a^{k-1}-|N_{H_k(U_0,T)}(v)\setminus N_G(v)|\ge a^{k-1}-(2k)^{k+4}\eta N^{k-1}.
	\end{aligned}
	\]
	Since $N\le n\le2ka$ and $\eta<(2k)^{-2k-4}$,
	$$\delta_1(J_0)
	\ge\bigl(1-(2k)^{2k+3}\eta\bigr)a^{k-1}
	>\left(1-\frac1{2k}\right)a^{k-1}\geq(\frac{k-1}{k}a^{k-1}+ka^{k-2}).$$
	By Theorem~\ref{thm:hps}, there is a perfect matching $M_2$ of $J_0$. Thus $M_0\cup M_1\cup M_2$ is a matching of $H$ of size
	$|M_0|+|M_1|+|M_2|=(b+1)+|M_1|+a=s$.
\end{proof}

\section{Proof of Theorem \ref{thm:main}}\label{sec:extension}
In this section, we prove Theorem~\ref{thm:main}.
We first establish the result for 
$\frac{n-c_0}{k}\leq s\leq \frac{n}{k}$ and then extend $s$ to the full range
by adjoining universal vertices, where $c_0$ is a fixed constant. We begin by combining Lemmas~\ref{lem:localcover}
and~\ref{prop:nonextremal} under the
minimum degree condition.


\begin{lemma}\label{lem:largeminimum}
	Let $k,c$ be positive integers such that $k\ge3$ and let $\alpha$ be a constant such that $0<\alpha<1$.
	There exists $s_0=s_0(k,c,\alpha)$ such that for integers $s\ge s_0$ and $n=ks-1+c$ the following holds.
	If an $n$-vertex $k$-graph $H$ satisfies
	$\delta_1(H)\ge\alpha\binom{n-1}{k-1}$ and	$\sigma_k(H)>kD(n,k,s)$,
	then $\nu(H)\ge s$.
\end{lemma}
\begin{proof}
	Let $\eta>0$ be given by Lemma~\ref{lem:localcover}
for $k$ and $c$. Apply Lemma~\ref{prop:nonextremal}
with parameters $k,\alpha/2,\eta$, and let $\delta>0$
be the resulting constant.
Choose $s_0$ sufficiently large in terms of $k,c,\alpha$.
	
	Let $Z\subseteq V(H)$ be a set of size $c-1$ and $J=H-Z$. Write $N=|V(J)|=ks$.
	If $J$ is $\eta$-extremal, then there is a matching of size $s$ in $H$ by Lemma~\ref{lem:localcover}. We may therefore assume that $J$ is not $\eta$-extremal.
	Since $n$ is sufficiently large, we have
	\[
	\begin{aligned}
		\delta_1(J)\ge\delta_1(H)-(c-1)\binom{n-2}{k-2}\ge\left(\alpha-\frac{(c-1)(k-1)}{n-1}\right)
		\binom{n-1}{k-1}\ge\frac{\alpha}{2}\binom{n-1}{k-1}
		\ge\frac{\alpha}{2}\binom{N-1}{k-1}.
	\end{aligned}
	\]
	Since every nonedge of $J$ is a nonedge of $H$,
\[
	\begin{aligned}
		\sigma_k(J)
		\ge\sigma_k(H)
		-k\left(\binom{n-1}{k-1}-\binom{N-1}{k-1}\right)\ge kD(N,k,N/k)-\delta\binom{N-1}{k-1}.
	\end{aligned}
	\]
    By Lemma~\ref{prop:nonextremal}, $J$ has a perfect
matching of size $s$, which is also a matching in $H$.
\end{proof}

\begin{theorem}\label{thm:fixedsurplus}
		Let $k,c$ be integers such that $k\ge3$ and $c\ge1$. There exists
		$s_0=s_0(k,c)$ such that for integers $s\ge s_0$ and $n=ks-1+c$ the following holds. Let $H$ be an $n$-vertex
		$k$-graph. If
		\[
		\sigma_k(H)>\max\left\{kD(n,k,s),\ (k-c)\binom{ks-2}{k-1}\right\},
		\]
		then $\nu(H)\ge s$.
	\end{theorem}
\begin{proof}
	Let $\alpha=(128k^3)^{-(k-1)}$. Choose $s_0$ sufficiently large for all subsequent
applications of Lemma~\ref{lem:largeminimum}
and the estimates below.
	By Lemma~\ref{lem:largeminimum}, we may assume that $\delta_1(H)<\alpha\binom{n-1}{k-1}$.
	Let $v\in V(H)$ such that
	$d_H(v)<\alpha\binom{n-1}{k-1}$, and define
	$
	B=\left\{u\in V(H):
	d_H(u)<D(n,k,s)+\frac1{8k}\binom{n-1}{k-1}\right\}$.
	Thus $v\in B$.
	Since $k,c$ are fixed and $s=(n-c+1)/k$,
for sufficiently large $n$ we have
    \begin{equation}\label{eq:fixedcoverlower}
		\begin{aligned}
			D(n,k,s)=\binom{n-1}{k-1}-\binom{n-s}{k-1}
			>(\frac{1}{2}-\frac{1}{4k})\binom{n-1}{k-1}.
		\end{aligned}
	\end{equation}
	For every $A\in\binom{B\setminus\{v\}}{k-1}$, we have
	\[
	\begin{aligned}
		d_H(v)+\sum_{u\in A}d_H(u)
		<(k-1)D(n,k,s)
		+\left(\alpha+\frac{k-1}{8k}\right)\binom{n-1}{k-1}
		<kD(n,k,s)<\sigma_k(H).
	\end{aligned}
	\]
	Thus $A\cup\{v\}$ is an edge. If $|B|>n/(64k^3)$, then
	\[
	\begin{aligned}
		\alpha>\frac{d_H(v)}{\binom{n-1}{k-1}}
		\ge\frac{\binom{|B|-1}{k-1}}{\binom{n-1}{k-1}}=\prod_{i=1}^{k-1}\frac{|B|-i}{n-i}\ge\left(\frac{n/(64k^3)-(k-1)}{n}\right)^{k-1}\ge\left(\frac1{128k^3}\right)^{k-1}=\alpha
	\end{aligned}
	\]
	for sufficiently large $n$, a contradiction. Hence
		$|B|\le\frac{n}{64k^3}$.

	Let
	\[
	L=\left\{u\in V(H):
	d_H(u)<\frac1{8k}\binom{n-1}{k-1}\right\}.
	\]
	Then $v\in L\subseteq B$. 
	
	\smallskip\noindent
	\textbf{Claim 5.} \emph{There is a matching $M$ in $H$ such that $|L\setminus V(M)|\leq c-1$ and every edge of $M$ intersects $L$.}
    
	\smallskip\noindent
	\emph{Proof of Claim 5.}
	Let $e$ be a $k$-set containing at
	least $k-1$ vertices of $L$. By \eqref{eq:fixedcoverlower},
	\[
	\begin{aligned}
		\sum_{u\in e}d_H(u)\le(1+\frac{k-1}{8k})\binom{n-1}{k-1}<k(\frac{1}{2}-\frac{1}{4k})\binom{n-1}{k-1}
		<kD(n,k,s)<\sigma_k(H).
	\end{aligned}
	\]
	Thus every $k$-set $e$ containing at least $k-1$ vertices of $L$ is an edge. Since $|L|\leq |B|\leq \frac{n}{64k^3}$, there exists a matching $M_1$ of size $\lfloor\frac{|L|}{k-1}\rfloor$ such that $|e\cap L|=k-1$ for each $e\in M_1$. Let $R=L\setminus V(M_1)$. Thus $|R|\le k-2$. If $|R|\leq c-1$, let $M=M_1$ and the proof of Claim 5 is done. 
	So we assume that $k-2\geq |R|\ge c$ in the following proof. 
	Let $f$ be a $k$-set containing at least $c+1$ vertices of $L$. Since
	\[
	\begin{aligned}
		\sigma_k(H)
		>(k-c)\binom{ks-2}{k-1}
		>(k-c-1+\frac{c+1}{8k})\binom{n-1}{k-1}\geq \sum_{u\in f}d_H(u),
	\end{aligned}
	\]
	 every $k$-set containing at least $c+1$ vertices
	of $L$ is an edge.
	If $|R|\ge c+1$, choose
$Y\in\binom{V(H)\setminus (B\cup V(M_1)}{k-|R|}$ and take
$M=M_1\cup\{R\cup Y\}$.
	
	It remains to consider the case $|R|=c$. If $M_1\neq\varnothing$, choose $e\in M_1$ and $r\in R$.
Since $|e\cap L|=k-1$, the set
    $e^*:=(e\cap L)\cup\{r\}$
is a $k$-set contained in $L$, and hence is an edge of $H$.
As $r\notin V(M_1)$, 
$    M:=(M_1\setminus\{e\})\cup\{e^*\}$
is a matching. Every edge of $M$ intersects $L$, and
    $|L\setminus V(M)|=|R|-1=c-1$.
This proves the claim in this case.
	We may therefore assume that $M_1=\varnothing$. Then $L=R$ and $|L|=c$.
	If every vertex of $L$ were isolated, a $k$-set containing $L$ would then
	be a nonedge with degree sum at most
	$(k-c)\binom{n-c-1}{k-1}=(k-c)\binom{ks-2}{k-1}<\sigma_k(H)$,
	a contradiction. Thus some edge $e'$ intersects $L$.
Taking $M=\{e'\}$ completes the proof of the claim. \claimqed
	\smallskip

	We next extend $M$ greedily to cover $B\setminus L$.
By Claim 5, every edge of $M$ intersects $L\subseteq B$
and $|L\setminus V(M)|\le c-1$.
Suppose that the current matching does not cover
$B\setminus L$, and let
$u\in(B\setminus L)\setminus V(M)$.
Since every edge of the current matching intersects $B$,
we have $|M|\le |B|$.
As $u\notin L$, the number of edges containing $u$
and avoiding $V(M)$ is at least
\[
\begin{aligned}
d_H(u)-k|M|\binom{n-2}{k-2}
\ge \frac1{8k}\binom{n-1}{k-1}
 -k|B|\binom{n-2}{k-2}
>0.
\end{aligned}
\]
Hence there exists an edge containing $u$ and
disjoint from $V(M)$, which can be added to $M$.
Repeating this process, we obtain a matching covering
$B\setminus L$. We still denote this matching by $M$.
Since every edge of $M$ intersects $B$, we have
$|M|\le |B|$.
Furthermore, as the matching is only enlarged,
the inequality $|L\setminus V(M)|\le c-1$ remains valid.
Consequently,
$B\setminus V(M)\subseteq L\setminus V(M)$ and
$|B\setminus V(M)|\le c-1$.
	
	Let $D_0\subseteq V(H)\setminus V(M)$ be a set of size $(c-1)$ such that
	$B\setminus V(M)\subseteq D_0$ and let $H_0=H-(V(M)\cup D_0)$.
	Then
	\[
	\begin{aligned}
		N_0:=|V(H_0)|=n-k|M|-(c-1)\ge n-k|B|-(c-1)
		\ge\left(1-\frac1{32k^2}\right)n\ge\frac n2
	\end{aligned}
	\]
	and
	\[
	\begin{aligned}
		D(N_0,k,N_0/k)
		=\sum_{i=1}^{N_0/k-1}\binom{N_0-i-1}{k-2}\le\sum_{i=1}^{s-1}\binom{n-i-1}{k-2}
		=D(n,k,s).
	\end{aligned}
	\]
	Note that all vertices of $H_0$ lie outside $B$. Hence
	\begin{equation}\label{eq:fixedresidualdegree}
		\begin{aligned}
			\delta_1(H_0)
			&\ge D(n,k,s)+\frac1{8k}\binom{n-1}{k-1}
			-(k|B|+c-1)\binom{n-2}{k-2}\\
			&>D(n,k,s)+\frac1{16k}\binom{n-1}{k-1}\\
			&\ge D(N_0,k,N_0/k)+\frac1{16k}\binom{N_0-1}{k-1}.
		\end{aligned}
	\end{equation}
	 By \eqref{eq:fixedcoverlower} and
	\eqref{eq:fixedresidualdegree}, 
	$\delta_1(H_0)>\frac14\binom{n-1}{k-1}
	\ge\alpha\binom{N_0-1}{k-1}$ and
	$\sigma_k(H_0)\ge k\delta_1(H_0)>kD(N_0,k,N_0/k)$.
	Applying Lemma~\ref{lem:largeminimum} to $H_0$, we obtain a perfect matching of $H_0$.
Together with $M$, this yields a matching of size $s$ in $H$.
\end{proof}

For integers $k\ge3$ and $s\ge1$, define
\[
c(s,k)=\min\left\{
d\in\{0,1,\ldots,k-1\}:
kD(ks+d,k,s)\ge
(k-d-1)\binom{ks-2}{k-1}
\right\}.
\]
This minimum exists since the inequality holds for $d=k-1$.
For fixed $k,s$, the left-hand side is nondecreasing in $d$,
whereas the right-hand side is nonincreasing. Thus the
inequality holds for every integer $d\ge c(s,k)$.

\begin{theorem}\label{thm:fullcover}
	Let $k$ be an integer such that $k\ge3$. There exists $n_0=n_0(k)$
	such that the following holds for all integers $n,s$
	with $n\ge n_0$, $s\ge1$ and $n\ge ks+c(s,k)$.
	If $H$ is an $n$-vertex $k$-graph satisfying
	$\sigma_k(H)>kD(n,k,s)$,
	then $H$ contains a matching of size $s$.
\end{theorem}
\begin{proof}
	If $n\le ks+k-2$, let $c'=n-ks+1$. Then $1\le c'\le k-1$.
	By the definition of $c(s,k)$ and monotonicity,
	\[
	\sigma_k(H)>kD(n,k,s)
	=\max\left\{kD(n,k,s),\,
	(k-c')\binom{ks-2}{k-1}\right\}.
	\]
	Since there are finitely many values of $c'$, by choosing $n_0$ sufficiently
	large in terms of $k$, Theorem~\ref{thm:fixedsurplus} yields a matching of size $s$
	in $H$. Hence assume that $n\ge ks+k-1$ in the following proof.
	
	Choose the unique integer $c'\equiv n-ks+1\pmod{k-1}$ such that
	\begin{equation}\label{eq:residues}
		c'-(k-1)\in\{1,2,\ldots,k-1\}.
	\end{equation}
	Since $n-ks+1\ge k$, and $c'$ is the least integer above
	$k-1$ in its residue class,
	$a:=\frac{n-ks+1-c'}{k-1}$
	is a nonnegative integer.
	
	Let $A$ be a set of $a$ new vertices disjoint from $V(H)$,
	and define $H^+$ on $V(H)\cup A$ by
	\[
	E(H^+)=E(H)\cup
	\left\{e\in\binom{V(H)\cup A}{k}:e\cap A\ne\varnothing\right\}.
	\]
	Let $N=n+a$ and $S=s+a$. Then
	$N=kS-1+c'$ and $N-S=n-s$.
	The non-edges of $H^+$ are precisely those of $H$, and every
	vertex $v\in V(H)$ satisfies
	$d_{H^+}(v)
	=d_H(v)+\binom{N-1}{k-1}-\binom{n-1}{k-1}$.
	Consequently,
	\begin{equation}\label{eq:extensionexcess}
		\begin{aligned}
			\sigma_k(H^+)-kD(N,k,S)
			&=\sigma_k(H)
			+k\left(\binom{N-1}{k-1}-\binom{n-1}{k-1}\right)-k\left(\binom{N-1}{k-1}-\binom{N-S}{k-1}\right)\\
			&=\sigma_k(H)-kD(n,k,s)>0.
		\end{aligned}
	\end{equation}
	
	There are only $k-1$ possible values of $c'$ in
	\eqref{eq:residues}, and each satisfies $k\le c'<2k$.
	Hence
	\[
	(k-c')\binom{kS-2}{k-1}\le0\le kD(N,k,S).
	\]
	Moreover, $S=\frac{N-c'+1}{k}>\frac nk-2$.
	Thus, by choosing $n_0$ sufficiently large in terms of $k$
	alone, Theorem~\ref{thm:fixedsurplus} applies for every
	possible $c'$.
	Together with \eqref{eq:extensionexcess}, we obtain
	\[
	\sigma_k(H^+)>
	\max\left\{kD(N,k,S),\,
	(k-c')\binom{kS-2}{k-1}\right\}.
	\]
	Theorem~\ref{thm:fixedsurplus} therefore gives a matching
	$M^+$ of size $S$ in $H^+$.
	Since the edges of $M^+$ are pairwise disjoint, at most
	$|A|=a$ of them meet $A$. Every remaining edge belongs to
	$H$. Deleting the edges that intersect $A$ therefore leaves at
	least $S-a=s$ pairwise disjoint edges in $H$, as required.
\end{proof}

\begin{proof}[Proof of Theorem~\ref{thm:main}]
If $n\ge ks+c(s,k)$, then $\sigma_k(H)>kD(n,k,s)$.
By Theorem~\ref{thm:fullcover}, $\nu(H)\geq s$. Hence we may assume that $n<ks+c(s,k)$. Since $c(s,k)\le k-1$,
we have
$1\le c\le k$.
Thus $c$ takes only finitely many values depending on $k$. We may therefore choose $n_0(k)$
sufficiently large so that, whenever $n\ge n_0(k)$,
$s$ satisfies the requirement of
Theorem~\ref{thm:fixedsurplus} for every $1\le c\le k$.
By
Theorem~\ref{thm:fixedsurplus}, there is a matching of size $s$ in $H$.
\end{proof}

\section*{Declaration of the use of AI}
ChatGPT was used to assist in developing the proof of Theorem \ref{thm:alternative}, as well as preparing the manuscript. All AI-suggested arguments were independently verified, revised and streamlined by the authors, who take full responsibility for the correctness of all mathematical results presented in this paper.

\end{document}